\documentclass[a4paper]{article}

\usepackage{amsfonts}
\usepackage{amsmath}
\usepackage{amsthm}\usepackage{amssymb}
\usepackage[mathscr]{eucal}
\usepackage{graphics}
\usepackage[cp866]{inputenc}
\usepackage[english]{babel}
\usepackage{makeidx}
\usepackage{latexsym,amsfonts,amssymb,amsmath,longtable,amsthm}
\input amssym.def
\usepackage{latexsym}
\usepackage{graphicx}
\usepackage[all]{xy}
\usepackage{amsfonts}
\usepackage{amsmath}
\usepackage{mathrsfs}
\usepackage{color}
\usepackage{ulem}
\usepackage[table]{xcolor}

\newtheorem{ttt}{Theorem}

\newtheorem{lem}{Lemma}
\newtheorem{cor}{Corollary}

 \def\Xint#1{\mathchoice
   {\XXint\displaystyle\textstyle{#1}}%
   {\XXint\textstyle\scriptstyle{#1}}%
   {\XXint\scriptstyle\scriptscriptstyle{#1}}%
   {\XXint\scriptscriptstyle\scriptscriptstyle{#1}}%
   \!\int}
\def\XXint#1#2#3{{\setbox0=\hbox{$#1{#2#3}{\int}$}
     \vcenter{\hbox{$#2#3$}}\kern-.5\wd0}}

\def\dashint{\Xint-}

 \usepackage[pagebackref,
            colorlinks,linkcolor={blue},citecolor={red}]
           {hyperref}

\def\0{{\mathbf{0}}}

\def\R{{\mathbb R}}
\def\th{{\theta}}

\def\H{{\mathcal{H}}}

\def\N{{\mathbb{N}}}

\def\ue{{\mathbf{u}}}
\def\loc{{\mathrm{loc}}}

\def\div{\hbox{\rm div}\,}
\def\arg{{\mathrm{arg}\,}}

 \def\into{\int\limits_}

\def\u{\mathbf u}
\def\brf{\bar f}
\def\bu{\bar{\mathbf u}}

\def\e{\varepsilon}
\def\a{\mathbf a}

\def\div{\hbox{\rm div}\,}

\newcommand{\diam}{\mathop{\mathrm{diam}}}

\def\et{{\mathbf e_{\tilde\theta}}}

\def\tl{{\tilde\theta}}
\def\tvv{{\tilde\varphi}}

\begin{document}
\title{Leray's plane steady state solutions are nontrivial}

\author{ Mikhail  Korobkov,  Konstantin Pileckas, 
and Remigio Russo}

\maketitle

\medskip

\noindent  
{\bf Affiliations:}  School of Mathematical Sciences,
Fudan University, Shanghai 200433, China; and Voronezh State University,
Voronezh 394018, Russia;\\
e-mail: {\it korob@math.nsc.ru}

\medskip

\noindent Institute of Applied Mathematics, Vilnius University, Naugarduko Str.,
24, Vilnius, 03225  Lithuania; \\
e-mail: {\it  konstantinas.pileckas@mif.vu.lt}

\medskip

\noindent Dipartimento di Matematica e Fisica
Universit\`a degli studi della Campania "Luigi Vanvitelli," viale
Lincoln 5, 81100, Caserta, Italy;\\
e-mail: {\it
remigio.russo@unicampania.it}

\medskip

\noindent{\bf Acknowledgment.}
 M. Korobkov was partially supported by the Ministry of Education
and Science of the Russian Federation (Grant 14.Z50.31.0037).

 The research of K. Pileckas was funded by the grant No. S-MIP-17-68 from
the Research Council of Lithuania.

\medskip


\begin{abstract}
 We study solutions to the obstacle problem for the
stationary Navier--Stokes system in a~two dimensional exterior
domain (flow past a prescribed body). We prove that the classical Leray solution to this problem is always nontrivial. No additional condition (on symmetry or smallness, etc.)  is assumed.
 This is  a complete extension of a~classical result  of  C.J. Amick  (Acta Math. 1988) where   nontriviality was  proved under symmetry assumption.
 
\noindent \medskip


\medskip

\end{abstract}

\setcounter{section}{0}

\section{Introduction}

\setcounter{equation}{0}
\setcounter{ttt}{0}
\setcounter{lem}{0}

 Let $\Omega$ be an exterior domain in $\R^2$ with compact boundary~$\partial\Omega=\bigcup\limits_{i=1}^N\Gamma_i$, where
$\Gamma_i$ are smooth disjoint curves, homeomorphic to the~circle. In particular, $\Omega\supset\R^2\setminus B$, where $B$ is the disk of radius $R_0$ centered at the origin with
$\partial \Omega\subset B$.

One of the most difficult and still  open problem in the
theory of  the stationary Navier--Stokes equations, initiated  by J. Leray in the famous paper
of 1933~\cite{Ler}, concerns the existence of a solution to the  {\it flow around an~obstacle\/}
(see also~\cite{Galdibook}):
\begin{equation}
\label{SNS-REM}
\left\{\begin{array}{r@{}l}
-\nu \Delta{\u}+(\u\cdot\nabla)\u+\nabla p  & {} ={\bf 0}\qquad \hbox{\rm in } \Omega, \\[2pt]
\div{\u} & {} =0\,\qquad \hbox{\rm in } \Omega\\[2pt]
\u& {} ={\bf 0}\qquad \hbox{\rm on } \partial\Omega,\\[2pt]
\u(z) & {} \to \u_0\quad\mbox{ as \ }|z|\to\infty,
 \end{array}\right.
\end{equation}
where $\u$ and $p$ are the unknown velocity  and pressure fields,  $\nu$ denotes the kinematical viscosity coefficient, and $\u_0\in\R^2$ is
a~nonzero constant vector (prescribed velocity  at~infinity).

  Leray  suggested~\cite{Ler}  the following  elegant  approach to this problem which was called method of  ``invading domains''.
Denoting by $\u_k$ the solution to   the problem
\begin{equation}
\label{SNS-REM3}
\left\{\begin{array}{r@{}l}
-\nu \Delta{\u_k}+(\u_k\cdot\nabla)\u_k+ \nabla p_k  & {} ={\bf 0}\qquad \hbox{\rm in } \Omega_{k}, \\[2pt]
\div{\u_k} & {} =0\,\qquad \hbox{\rm in } \Omega_{k},  \\[2pt]
{\u_k} & {} = \mathbf0\,\qquad \hbox{\rm on } \partial\Omega,  \\[2pt]
\u_k & {}= \u_0\;\;\quad\mbox{for \ }|z|=R_k.
 \end{array}\right.
\end{equation}
on the  intersection $\Omega_k$ of $\Omega$ with the disk $B_{R_k}$ of   radius $R_k\ge k(\gg R_0)$,  whose existence he  proved before, Leray showed that
the sequence $\u_k$ satisfies the   estimate
$\int_\Omega|\nabla\u_k|^2\le c$ for some positive constant $c$ independent of $k$. Hence, he observed that it is possible to extract a subsequence $\ue_{k_n}$ which  weakly  converges  to a solution  $\ue_L$ of problem~$(\ref{SNS-REM})_{1,2,3}$ with  $\int_\Omega|\nabla\u_L|^2<+\infty$.
This solution   was later called {\it Leray's solution\/}  (see, $e.g.$, \cite{Amick}\,).

An arbitrary  solution $\u$ to the Navier--Stokes equations
\begin{equation}
\label{SNS}
\left\{\begin{array}{r@{}l}
-\nu \Delta{\u}+(\u\cdot\nabla)\u+\nabla p  & {} ={\bf 0}\qquad \hbox{\rm in } \Omega, \\[2pt]
\div{\u} & {} =0\,\qquad \hbox{\rm in } \Omega.
 \end{array}\right.
\end{equation}
having  the~finite Dirichlet integral
\begin{equation}
\label{SxxS}
\into\Omega|\nabla\u|^2<+\infty,
\end{equation}
is called today {\it D--solution\/} \cite{Galdibook}.   As is well known ($e.g.$, \cite{OAL}), such solutions are
real--analytic in $\Omega$.

 As far as condition(\ref{SNS-REM})$_4$ is concerned, Leray limited himself to observe that, while in three dimensional problem~(\ref{SxxS}) it is sufficient to guarantee the attainability of the limit $	\u_0$ at infinity (at least in a mean square sense) as a consequence of the inequality $\|r^{-1}(\u-\u_0)\|_{L^2(\Omega)}\le 4\|\nabla\u\|_{L^2(\Omega)}$, in the two dimension case the corresponding inequality $\| (r\log r)^{-1}(\u-\u_0)\|_{L^2(\Omega)}\le c\|\nabla\u\|_{L^2(\Omega)}$ does not imply any type of convergence. Leray concluded that one should not be surprised of this phenomenon, in view of the~{\it Stokes paradox},  $i.e.$,   the   system obtained from  (\ref{SNS-REM}) removing the nonlinear term, namely,
 \begin{equation}
\label{Stok}
\left\{\begin{array}{r@{}l}
-\nu \Delta{\u}+\nabla p  & {} ={\bf 0}\qquad \hbox{\rm in } \Omega, \\[2pt]
\div{\u} & {} =0\,\qquad \hbox{\rm in } \Omega\\[2pt]
\u& {} ={\bf 0}\qquad \hbox{\rm on } \partial\Omega,\\[2pt]
\u(z) & {} \to \u_0\quad\mbox{ as \ }|z|\to\infty,
 \end{array}\right.
\end{equation}
 does not admit  a~solution~(see, $e.g.$, \cite{RR10}).

 The problem  of the asymptotic behaviour at infinity of Leray's solution $(\u_L, p_L)$ was tacked by D.~Gilbarg \& H. Weinberger  in 1974 \cite{GW1}. They proved that
 $\u_L$ is bounded, there are a scalar $p_0$ and a constant vector $\u_\infty$ such that
\begin{equation}\label{prc}
\lim_{|z|\to+\infty}p_L(z)= p_0
\end{equation}
(one can choose, say, $p_0=0$\,),
\begin{equation}
\label{GGWW}
\displaystyle\lim_{|z|\to+\infty}\into0^{2\pi}|\u_L( r,\theta)-\u_\infty|^2d\theta=0,
\end{equation}
and
\begin{equation}
\label{SzzSL}
\begin{array}{ l}
 \omega(z)  =o(r^{-3/4}),\\[2pt]
\int\limits_{\Omega}r|\nabla\omega(z)|^2<\infty,
 \end{array}
\end{equation}
where $r=|z|$ and
$$
\omega=\partial_2u_{L1}-\partial_1u _{L2}
$$
is the vorticity.

In 1988 C.J. Amick \cite{Amick} proved that a $D$--solution to the  problem of a  flow around an obstacle  (\ref{SNS-REM})$_{1.2.3}$  has the following asymptotic  properties:
\begin{itemize}
\item[(i)] \, $\u$ is bounded and, as a consequence, it satisfies
(\ref{GGWW})--(\ref{SzzSL});
\item[(ii)] \,the
total head pressure $\Phi=p+\frac12|\u|^2$ and the absolute value
of the velocity $|\u|$ have the uniform limit at infinity, i.e.,
\begin{equation}
\label{dom6} |\u(r,\th)|\to|\u_\infty|\qquad\mbox{ as }\
r\to\infty,
\end{equation}
where $\u_\infty$ is the~constant vector from the
condition~(\ref{GGWW});
\item[(iii)] if $\partial\Omega$ is symmetric with respect to the $x_1$--axis, and $\u=(u_1,u_2)$ is also symmetric, i.e., if $u_1$ is even and $u_2$ is odd with respect to $x_1$, then $\u$ converges uniformly at infinity to a constant vector $\mu\mathbf e_1$, for some scalar~$\mu$. Moreover, the Leray procedure yields a {\it nontrivial } (i.e.,  not identically zero)  
symmetric  solution.
\end{itemize}

In the present paper we prove that in general case (without any additional symmetry assumptions) the Leray solution to the the problem of the flow around obstacle is
{\it always} nontrivial.

\begin{ttt}
\label{TT2} {\sl  Let $\Omega$ be an exterior domain in $\R^2$ with smooth compact boundary, $\nu>0$ and ${\bf 0}\ne\u_0\in\R^2$.
Take a sequence $\u_k$ of solutions to system~(\ref{SNS-REM3}), and take further arbitrary weakly convergent subsequence $\u_{k_n}\rightharpoonup \u$. Then the limiting solution~$\u$ to~$($\ref{SNS-REM}$)_{1,2,3}$ is nontrivial (i.e., $\u$ is not identically zero\,).   In particular the Leray solution is nontrivial.}
 \end{ttt}

  Moreover, we proved a kind of complementary result.
\begin{ttt}
\label{TT3} {\sl  Let $\Omega$ be an exterior domain in $\R^2$ with smooth compact boundary, $\nu>0$ and let $\a\in\R^2$ be a nonzero constant vector.
Take a sequence $\u_k$ of solutions to the~system
\begin{equation}
\label{SNS-S}
\left\{\begin{array}{r@{}l}
-\nu \Delta{\u_k}+(\u_k\cdot\nabla)\u_k+ \nabla p_k  & {} ={\bf 0}\qquad \hbox{\rm in } \Omega_{k}, \\[2pt]
\div{\u_k} & {} =0\,\qquad \hbox{\rm in } \Omega_{k},  \\[2pt]
{\u_k} & {} =\a\,\qquad \hbox{\rm on } \partial\Omega,  \\[2pt]
\u_k & {}= {\bf 0}\;\qquad\mbox{for \ }|z|=R_k,
 \end{array}\right.
\end{equation}
and take further arbitrary weakly convergent subsequence $\u_{k_n}\rightharpoonup \u$. Then the~limiting solution~$\u$  is nontrivial.
In other words, the solution to the system
$$
\left\{\begin{array}{r@{}l}
-\nu \Delta{\u}+(\u\cdot\nabla)\u+\nabla p  & {} ={\bf 0}\qquad \hbox{\rm in } \Omega, \\[2pt]
\div{\u} & {} =0\,\qquad \hbox{\rm in } \Omega\\[2pt]
\u& {} =\a\,\qquad \hbox{\rm on } \partial\Omega,
\end{array}\right.
$$
obtained by the Leray method is nontrivial, $i.e.$, $\u\ne\a$.}

 \end{ttt}

  Note, that for linear case ($e.g.$, for Stokes system~(\ref{Stok})\,) the~assertions of Theorems~\ref{TT2}--\ref{TT3} are evidently equivalent. But of course it does not hold in general for nonlinear systems.

 \medskip

 Recently \cite{kpr-arxiv}, \cite{kpr-18}   we  proved the following result for general~$D$-solutions.

\begin{ttt}[\cite{kpr-18}]
\label{T2} {\sl Let $\u$ be a $D$-solution to the Navier--Stokes
system~(\ref{SNS})
in the exterior domain~$\Omega\subset\R^2$.
Then $\u$ converges uniformly at infinity, i.e.,
\begin{equation}
\label{as-0}
\u(z)\to \u_\infty\qquad\mbox{ uniformly as \ }|z|\to\infty,
\end{equation}
where $\u_\infty\in\R^2$ is some constant vector. }
 \end{ttt}

  By virtue of Theorem \ref{T2}, for every Leray solution there is a constant vector
$\ue_\infty\in\R^2$ such that \eqref{as-0} holds.  However,
 the~desired equality~$\ue_\infty=\ue_0$ is still~an open question. We even do not know whether $\ue_\infty$ is nonzero if $\ue_0\neq {\bf 0}$.  So, the problem of the flow around an~obstacle  remains   still open.

The same open problem exists in~Theorem~\ref{TT3}: we know that the corresponding Leray solution converges uniformly to some constant vector
$\ue_\infty\in\R^2$, but   we were not able to prove the expected equality that  $\ue_\infty={\bf 0}$.

 \medskip

{\bf Some additional historical remarks.} \ Note, that thirty years after Leray,
H.~Fujita~\cite{Fujita} by means of  different techniques proved the  existence of a $D$--solution to (\ref{SNS-REM})$_{1.2.3}$. Due of a lack of a uniqueness theorem, the Leray and Fujita  solutions are not comparable.

Recall also  the~amazing  discovery of R. Finn and D.R. Smith in~1967 \cite{FS} of the existence of a~solution to (\ref{SNS-REM})\, for $\nu$ sufficiently large (or, equivalently, for $\u_0$ sufficiently small). Their    approach  is completely different from  that of Leray. Nevertheless, their method does not allow to prove the existence of the solution to the~problem
\begin{equation}
\label{SNS-00}
\left\{\begin{array}{r@{}l}
-\nu \Delta{\u}+(\u\cdot\nabla)\u+\nabla p  & {} ={\bf 0}\qquad \hbox{\rm in } \Omega, \\[2pt]
\div{\u} & {} =0\,\qquad \hbox{\rm in } \Omega\\[2pt]
\u& {} =\a\,\qquad \hbox{\rm on } \partial\Omega,\\[2pt]
\u(z) & {} \to {\bf 0}\;\;\quad\mbox{ as \ }|z|\to\infty,
 \end{array}\right.
\end{equation}
for a~constant vector~$\a\in\R^2$. This problem is quite open {\bf even for small }vectors~$\a\ne{\bf 0}$ \,(the~existence was proved  for some nonconstant~$\a$
 under the assumption of symmetry with respect to both coordinate axes~\cite{KRMann}\,). Even the reduced problem~(\ref{SNS-00}${}_{1-3}$) is open for general boundary value~$\a$ ~(see, $e.g.$, \cite{AR} for the case of small fluxes).
More detailed survey of results concerning boundary value problems for stationary NS-system in plane exterior domains see, $e.g.$, in ~\cite{Galdibook} or in our recent papers~\cite{kpr-arxiv},  \cite{kpr-18}, \cite{kprplane}.

 Finally, let us describe shortly the main steps of the proof of Theorem~\ref{TT2} (for defineteness,
 take~$\u_0=(1,0)$ \,and \,$\nu=1$\,). The main ideas are rather simple.
 By Amick criterion~\cite{Amick}, the corresponding limiting Leray solution $\ue$ is trivial if and only if the convergence
 \begin{equation}
\label{SNS-D1}
\int\limits_{\Omega_k}|\nabla\u_k|^2\,dx\to0\qquad\mbox{ as }\ k\to\infty
\end{equation}
holds for the sequence of solutions~$\u_k$ to~(\ref{SNS-REM3}) in bounded domains~$\Omega_k=\Omega\cap B_{R_k}$. Suppose (\ref{SNS-D1})  to be  fulfilled. Then the~functions $(\u_k,p_k)$ and all its derivatives go to~$0$ uniformly on every bounded set; moreover, from~\cite{GW1} it follows that~$\sup\limits_{z\in\Omega'_k}|p_k|\to 0$, where $\Omega'_k=\Omega\cap B_{\frac34 R_k}$. It is well known, that
the Bernoulli pressure $\Phi_k=p_k+\frac12|\u_k|^2$ satisfies the maximum principle. Using these facts,
 it can be proved that the level lines ~$\Phi_k=t$ of the Bernoulli pressure are arranged as circles surrounding the origin; furthermore,
the~Bernoulli  pressure approximately equal to zero near~$\partial\Omega$, and it increases up  to~$\frac12$
near the large circle of a~radius~$R_k$.

The following  steps  are crucial in our arguments:

1) The~{\it direction} of the~velocity vector~$\u_k$ is under control of the Dirichlet integral
(it was proved in the~Gilbarg--Weinberger  paper~\cite{GW2}, see below lemma~\ref{LemmaGW} for the exact formulation of the result);

2) The vorticity $\omega_k(z)$ does not change sign between two level lines of the Bernoulli pressure~$\Phi_k$
(it  is proved using the~results of Amick~\cite{Amick}\,).

Using these important facts,  we prove that for the velocity~$\ue_k$ the following representation formulas hold:
\begin{equation}
\label{en20-int}
\bar\u_k(r)=|\bar\u_k(r)|\,(\cos\varphi_k(r),\sin\varphi_k(r)),
\end{equation}
where $\bar\u_k(r)$ means the mean value of~$\u_k$ over the circle of radius~$r$  and
\begin{equation}
\label{en21-int}
|\varphi_k(r)|\le\e_k\qquad\mbox{ for all sufficiently large } r\le R_k
\end{equation}
with $\e_k\to0$ as $k\to\infty$. Recall, that the gradient of the Bernoulli pressure satisfies the identity
\begin{equation}
\label{en24-int}
\nabla\Phi_k(z)=-\nabla^\bot\omega_k(z)+\omega_k(z)\cdot\ue_k^\bot(z),
\end{equation}
where $\u_k^\bot:=(-u_k^2,u_k^1)$.
The first term~in (\ref{en24-int}) is negligible (since the integral~$\int\limits_{\Omega'_k}r|\nabla\omega_k|^2$ is small, see~(\ref{SzzSL}${}_2$)\,). Using these
 facts we obtain the contradiction with the~geometrical structure of the level lines of the Bernoulli pressure~$\Phi_k$ described above and, as a consequence, the nontriviality of the Leray solution.

\section{Notations and preliminaries}
\setcounter{equation}{0}

By {\it a domain} we mean an open connected set. We use standard
notations for Sobolev  spaces \,$W^{k,q}(\Omega)$, where $k\in{\mathbb N}$,
$q\in[1,+\infty]$. In our notation we do not
distinguish function spaces for scalar and vector valued
functions; it is clear from the context whether we use scalar or
vector (or tensor) valued function spaces.

For $q\ge1$ denote by $D^{k,q}(\Omega)$ the set of functions $f\in
W^{k,q}_{\loc}(\Omega)$ such that
$\|f\|_{D^{k,q}(\Omega)}=\|\nabla^k f\|_{L^q(\Omega)}<\infty$.

We denote by $\H^k$ the
$k$-dimensional Hausdorff measure, i.e.,
$\H^k(F)=\lim\limits_{t\to 0+}\H^k_t(F)$,
where $$\H^1_t(F)=\bigl(\frac{\alpha_k}2\bigr)^k\inf\big\{\sum\limits_{i=1}^\infty \bigl({\rm
diam} F_i\bigr)^k:\, {\rm diam} F_i\leq t, F\subset
\bigcup\limits_{i=1}^\infty F_i\big\}$$
and $\alpha_k$ is a Lebesgue volume of the unit ball in~$\R^k$.

In particular, for a curve $S$ the value $\H^1(S)$ coincides with its length, and for sets $E\subset\R^2$ the $\H^2(E)$ coincides with the usual Lebesgue measure in~$\R^2$.

Also, for a curve $S$ by $\int\limits_Sf\,ds$ we denote the usual integral with respect to $1$-dimensional Hausdorff measure (=length). Further, for a set $E\subset\R^2$ by
$\int\limits_Ef(\xi)\,d\H^2$ we denote the integral with respect to the two-dimensional Lebesgue measure.  For convenience, we will write simply
$$\int\limits_E\!rf\qquad\mbox{ instead of }\qquad\int\limits_E|\xi|\, f(\xi)\,d\H_\xi^2.$$

Below we present some  results concerning the behavior of~$D$-functions.

 \begin{lem}
\label{lt1} {\sl Let $f\in D^{1,2}(\Omega)$ and assume that
$$\int\limits_{D}|\nabla f|^2<\e^2$$ for some $\e>0$ and for some ring
$D=\{z\in\R^2:r_1<|z-z_0|<r_2\,\}\subset\Omega$. Then the~estimate \begin{equation}
\label{est-1}
|\brf(r_2)-\brf(r_1)|\le \e\sqrt{\ln\frac{r_2}{r_1}}
\end{equation} holds, where $\brf$ is the mean value of~$f$ over the circle $S(z_0,r)$:
$$\brf(r):=\frac1{2\pi r}\int\limits_{|z-z_0|=r}f(z)\,ds.$$ }
\end{lem}

 \begin{lem}
\label{lt2} {\sl Fix a number $\beta\in(0,1)$. Let $f\in D^{1,2}(\Omega)$ and assume that
$$\int\limits_{D}|\nabla f|^2<\e^2$$ for some $\e>0$ and for  some ring
$D=\{z\in\R^2:\beta R<|z-z_0|<R\,\}\subset\Omega$. Then there exists  a number $r\in[\beta R, R]$ such that
 the
estimate \begin{equation}
\label{est-2}
\sup\limits_{|z-z_0|=r}|f(z)-\brf(r)|\le c_\beta\e
\end{equation} holds, where the constant~$c_\beta$ depends on $\beta$ only.   }
 \end{lem}

The proofs of above lemmas  are standard, see, $e.g.$, \cite{GW2} for the proofs of similar results.
Summarizing the results of these lemmas, we receive
\begin{lem}\label{lt3}
{\sl Under conditions of Lemma~\ref{lt2}, there exists $r\in [\beta R,R]$ such that
\begin{equation}
\label{est-3}
\sup\limits_{|z-z_0|=r}|f(z)-\brf(R)|\le \tilde c_\beta\e.
\end{equation}}
\end{lem}

The following result  was proved in~\cite[Theorem~4, page~399]{GW2}. It means, roughly speaking, that the~{\it direction} of the~velocity vector~$\u$ satisfying the Navier--Stokes system  is  controlled by  the Dirichlet integral.

\begin{lem}[\cite{GW2}]
\label{LemmaGW}
{\sl Let $\u$ be a $D$-solution to the Navier--Stokes
system~(\ref{SNS})
in the exterior domain~$\Omega\subset\R^2$, and let
$D=\{z\in\R^2:R_1<|z|<R_2\,\}\subset\Omega$ be some ring in  $\Omega$.
Denoted by $\bu(r)$ the mean value  of $\u$ over the circle $S_r$:
\begin{equation}
\label{eqv-c3}
\bu(r)=\frac1{2\pi r}\int\limits_{|\xi|=r}\u(\xi)\,ds
\end{equation}
and let $\varphi(r)$ be the argument of the complex number associated to the vector~$\bu(r)=(\bar u_1(r), \bar u_2(r))$, i.e., $\varphi(r)=\arg(\bar u_1(r)+i \bar u_2(r))$.
Assume also  that
$$|\bu(r)|\ge\sigma$$
for some positive constant~$\sigma>0$ and for all $r\in[R_1,R_2]$. Then  the following estimate
\begin{equation}
\label{eqv-c4--}
\sup\limits_{R_1< \rho_1\le\rho_2\le R_2}|\varphi(\rho_2)-\varphi(\rho_1)|\le \frac1{4\pi\sigma^2}\int\limits_{D}\biggl(\frac{1}{r}|\nabla\omega|+
|\nabla\u|^2\biggr)
\end{equation}
holds.}
\end{lem}

The following statement also follows from~\eqref{eqv-c4--}.

\begin{cor}
\label{LemmaGW-cor}
{\sl Under the~conditions of Lemma~\ref{LemmaGW}, the estimate
\begin{equation}
\label{ess-angle}
\sup\limits_{R_1< \rho_1\le\rho_2\le R_2}|\varphi(\rho_2)-\varphi(\rho_1)|
\le \frac1{2\sigma^2R_1}+
\frac1{4\pi\sigma^2}\int\limits_{D}\biggl( r\,|\nabla\omega|^2+
|\nabla\u|^2\biggr)
\end{equation}
holds.}
\end{cor}

\section{Proof of the main Theorem~\ref{TT2}.}
\setcounter{equation}{0}
\setcounter{lem}{0}

We prove Theorem~\ref{TT2} by getting a contradiction. Suppose the assumptions of Theorem~\ref{TT2} are fulfilled, but the statement is false, i.e.,
 there exists an increasing  sequence of radii $R_k\to+\infty$ and solutions $\u_k$ to the system~(\ref{SNS-REM3})
 such that $\u_{k_n}\rightharpoonup \u\equiv{\bf 0}$. By the result of Amick~\cite[Theorem~24, page~115]{Amick}, it is equivalent to the global convergence to zero of the Dirichlet integrals:
 \begin{equation}
\label{nont-1}
\int\limits_{\Omega_k}|\nabla\u_k|^2\to 0.
\end{equation}
 By
classical regularity results for $D$-solutions to the Navier--Stokes
system ($e.g.$, \cite{Galdibook}),  the functions~$\u_k$ and $p_k$ are $C^\infty$-smooth on the set~$\overline\Omega_k$ and real analytical inside~$\Omega_k$.  Moreover, (\ref{nont-1}) implies in particular, that for every compact set $E\subset\overline{\Omega}$
  \begin{equation}
\label{nont-2}
\sup\limits_{x\in E}|\nabla^j\u_k(x)|\to 0 \qquad\forall j=0,1,2,\dots\end{equation}
uniformly as $k\to\infty$, i.e., $\u_k$ and all its derivatives converges to zero as~$k\to\infty$ uniformly on every compact set.

Without loss of generality we may assume that
\begin{equation}
\label{nont-3}
\nu=1\qquad \mbox{ and }\qquad\u_0=(1,0).
\end{equation}

The proof consists of eight steps.

\medskip

{\sc Step 1}.  Denote $\Omega'_k=\Omega\cap B(0,\frac34 R_k)$. By results of \cite{GW1}--\cite{GW2} (see the proofs of Lemmas 2.2 and 3.2 in \cite{GW1}), the following estimate
$$
\int\limits_{\Omega_k\cap B(0, \frac{3}{4}R_k)} r |\nabla \omega_k|^2\leq c \int\limits_{\Omega_k} | \omega_k|^2\leq c \int\limits_{\Omega_k} |\nabla \u_k|^2
$$
holds with the constant $c$ independent of $k$. Hence the   assumption~(\ref{nont-1})   yields
 \begin{equation}\label{nont-5}
\int\limits_{\Omega'_k} r\,|\nabla\omega_k|^2\to 0.
\end{equation}
Moreover, from \cite{GW1}--\cite{GW2}  it follows that (see \cite[Lemma~2.4]{GW1}\,)
\begin{equation}\label{nont-7}
\sup\limits_{x\in\Omega'_k}|\ue_k(x)|\le C.
\end{equation}
and (see the proofs of Lemmas 2.3--2.6 in \cite{GW1} and, in particular, \cite[formula~(2.20)]{GW1}\,)
 \begin{equation}\label{nont-6}
\sup\limits_{x\in\Omega'_k}|p_k(x)|=\e_k.
\end{equation}
Here and everywhere below the equality  $a_k=\e_k$  means  that the~sequence~$a_k$ tends to $0$ as $k\to\infty$.

\

{\sc Step 2.} Denote
\begin{equation}\label{exp-1}
R_{0k}=\min\biggl\{r\ge R_0:|\bar\u_k(r)|=\biggl|\dashint_{S_r}\u_k\,ds\biggr|=\frac15\biggr\},
\end{equation}
where, as usual, $S_r=\{\xi\in\R^2:|\xi|=r\,\}$ is a circle.  From ~(\ref{nont-2}) we have, in particular,
\begin{equation}\label{nont-9}
R_{0k}\to+\infty\qquad\mbox{ as }\ k\to\infty.
\end{equation}
Moreover,  Lemma~\ref{lt1}   applied to $\u_k$  and the identity~$\bar\u_k(R_k)=(1,0)$  imply,  by virtue of \eqref{nont-1}, that
\begin{equation}\label{nont-10}
\frac{R_{k}}{R_{0k}}\to+\infty \qquad\mbox{ as }\ k\to\infty.
\end{equation}

{\sc Step 3.} By construction (see \eqref{nont-6}, \eqref{exp-1}),
\begin{equation}\label{nont-13}
\max\limits_{x\in S_{R_{0k}}}\Phi_k(x)>\frac12\biggl(\frac15\biggr)^2+\e_k=\frac1{50}+\e_k>\frac1{60}>\max\limits_{x\in \partial\Omega}\Phi_k(x),
\end{equation}
for sufficiently large~$k$, where, recall, $\Phi_k=p_k+\frac12|\u_k|^2$ is the Bernoulli pressure.
It is well known, that $\Phi_k$ satisfies the  identity
\begin{equation}\label{nont-11}
\Delta\Phi_k=\omega_k^2+\nabla\Phi_k\cdot\u_k,
\end{equation}
 and thus, it satisfies the classical strong maximum principle ($e.g.$, \cite{GW1}--\cite{GW2}\, or \cite{EvPDE}\,). From this property and ~(\ref{nont-13}) we obtain
\begin{equation}\label{nont-14}
\max\limits_{x\in S_{R}}\Phi_k(x)=\max\limits_{x\in \Omega_{R}}\Phi_k(x)\qquad\ \forall R\in[R_{0k},R_k],
\end{equation}
where, as usual, $S_R=\{x\in\R^2:|x|=R\,\}$  is a circle and $\Omega_R=\Omega\cap B_R=\{x\in\Omega:|x|<R\}$.
From  this maximum principle it follows that
\begin{equation}\label{nont-15}
\mbox{ \textsl{the function} \ \ }[R_{0k},R_k]\ni R\mapsto \max\limits_{x\in S_{R}}\Phi_k(x)\quad\mbox{\textsl{is strictly increasing}}.
\end{equation}
In particular, by (\ref{nont-13}),
\begin{equation}\label{nont-15-1}
 \max\limits_{x\in S_{R}}\Phi_k(x)>\frac1{60}\qquad \forall  R\in[R_{0k},R_k].
\end{equation}

{\sc Step 4.}  By results of  \cite{GW1}--\cite{GW2}(see the proof of Lemma 4.1 in \cite{GW1}, in particular, the proof  of the  estimate (4.5)\,) there exist a sequence of radii $R_{mk}$, $m=1,2,\dots,M=M(k)$, such that
\begin{eqnarray}
\label{nont-17}
2^{m-1}R_{0k}< R_{mk}<2^mR_{0k},\qquad  m=1,2,\dots, M;\\
\label{nont-18}
\frac13R_k<R_{Mk}<\frac23R_k;\\
\label{nont-19}
\sup\limits_{x\in S_{R_{mk}}}|\u_k(x)-\bar\u_k(R_{mk})|=\e_k.
\end{eqnarray}
Taking $k$ large enough and using  estimates~(\ref{nont-15-1}) for the Bernoulli pressure $\Phi_k$ and  estimates~(\ref{nont-6}) for the pressure~$p_k$, we conclude from~(\ref{nont-19}) that
\begin{equation}\label{nont-20}
\bigl|\bar\u_k(R_{mk})\bigr|=\biggl|\dashint_{S_{R_{mk}}}\u_k\,ds\biggr|>\sqrt{2\cdot \frac1{60}}+\e_k>\frac1{6}\qquad\forall m=1,2,\dots,M.
\end{equation}
Then, applying Lemma~\ref{lt1} and using the smallness of the Dirichlet integrals~(\ref{nont-1}) and the condition~(\ref{nont-20}), \eqref{nont-17}, we obtain
\begin{equation}\label{nont-21}
\bigl|\bar\u_k(R)\bigr|=\biggl|\dashint_{S_{R}}\u_k\,ds\biggr|>\frac1{7}\qquad\forall R\in[R_{0k},R_k].
\end{equation}
Denote by $\varphi_k(r)$ the angle direction of the vector~$\bar\u_k(r)$, i.e.,
\begin{equation}
\label{ess-angle-3}
\bar\u_k(r)=|\bar\u_k(r)|\bigl(\cos\varphi_k(r),\sin\varphi_k(r)\bigr).
\end{equation}
From Corollary~\ref{LemmaGW-cor} and relations \eqref{nont-1},
(\ref{nont-5}), (\ref{nont-18}), \eqref{nont-21} we have
\begin{equation}
\label{ess-angle-1}
\sup\limits_{R_{0k}< \rho_1\le\rho_2\le R_{Mk}}|\varphi_k(\rho_2)-\varphi_k(\rho_1)|=\e_k.\end{equation}
Further, Lemma~\ref{lt1} and the boundary conditions $\u_k(z)\big|_{|z|=R_k}= (1,0)$  yield
$$|\bar\u_k(r)-(1,0)|=\e_k\qquad\forall r\in[R_{Mk},R_k].$$
From the last two  formulas we conclude that
\begin{equation}
\label{ess-angle-2}
\sup\limits_{R_{0k}\le r\le R_k}|\varphi_k(r)|=\e_k.\end{equation}
Summarizing we can say that the formula
\eqref{ess-angle-3} holds with
\begin{equation}
\label{ess-angle-4}
|\bar\u_k(r)|\ge\frac17,\qquad|\varphi_k(r)|=\e_k\qquad\forall r\in[R_{0k},R_k].
\end{equation}

\

{\sc Step 5.} From the choice of the radii $R_{1k}$ and $R_{Mk}$ at Step~4  (see \eqref{nont-17}--\eqref{nont-18}), from estimates  \eqref{nont-6}, \eqref{nont-19}, and from the conditions
$|\bar\u_k(R_{0k})|=\frac15$  (see \eqref{exp-1}) and $|\bar\u_k(R_k)|=1$,
we have
\begin{equation}
\label{bpr1}
\sup\limits_{x\in S_{R_{1k}}}|\Phi_k(x)-\frac1{50}|\to0,
\end{equation}
\begin{equation}
\label{bpr2}
\sup\limits_{x\in S_{R_{Mk}}}\bigl|\Phi_k(x)-\frac1{2}|\to0
\end{equation}
as $k\to\infty$.  Without loss of generality we may assume that
\begin{equation}
\label{bpr3}
\Phi_k(x)<\frac1{45}\qquad\forall x\in{S_{R_{1k}}},
\end{equation}
\begin{equation}
\label{bpr4}
\Phi_k(x)>\frac1{3}\qquad\forall x\in{S_{R_{Mk}}}.
\end{equation}
Denote by $I$ the interval
$I=\bigl[\frac1{45},\frac13\bigr]$. By construction and by the classical Morse--Sard Theorem (which says  that the set of critical values of a $C^\infty$ function has zero Lebesgue measure) we  conclude  that
for almost all $t\in I$ the set
$$
\{x\in\R^2:R_{1k}\le|x|\le R_{Mk}, \ \ \Phi_k(x)=t\}
$$
is a finite disjoint union of smooth closed curves.  Moreover, every  of these curves is homeomorphic to the circle. (It follows from the fact, that the preimage of a non-critical value is a~smooth one dimensional manifold, and,  since $\Phi_k(x)\notin I$ for $x\in S_{R_{1k}}$  and $x\in S_{R_{Mk}}$,    this  manifold has no boundary.)   By evident topological reasons,
at least one of
these curves separate the circles $ S_{R_{1k}}$  and $S_{R_{Mk}}$. By maximum principle for the Bernoulli pressure $\Phi_k$ this  separating  curve is unique; denote it by $S_{k}(t)$.
In other words, we  have  proved  that
\begin{equation}
\label{lsa}
\begin{array}{c}
\mbox{\textsl{for almost all} $t\in I=\bigl[\frac1{45},\frac13\bigr]$ \textsl{there exists exactly one smooth curve}} \\
\mbox{ $S_k(t),$ \textsl{homeomorphic to the circle  separating} $ S_{R_{1k}}$ \textsl{from}  $S_{R_{Mk}}$, }\\
\mbox{ \textsl{and satisfying the identity}\;\;
$\Phi_k(x)\equiv t\quad\forall x\in S_{k}(t)$}. \end{array}
\end{equation}

\

{\sc Step 6.} Take numbers $t_1\in\bigl[\frac1{45},\frac1{40}\bigr]$, $t_2\in\bigl[\frac14,\frac13\bigr]$ which are regular values for all $\Phi_k$, \ $k=1,2,\dots$.
Then denote by $\Omega^s_k$ the bounded open subset of $\Omega_k$ satisfying
\begin{equation}
\label{bpr6}
\partial \Omega^s_k=S_k(t_1)\cup S_k(t_2).
\end{equation}
We claim that
\begin{equation}
\label{bpr8}
\mbox{\textsl{vorticity} $\omega_k(x)$ \textsl{does not change sign in} $\Omega^s_k$.}
\end{equation}
In order to prove this claim, consider the   auxiliary function
$$\gamma_k=\Phi_k-\omega_k\psi_k,$$
where $\psi_k$ is a stream function satisfying $\nabla\psi_k=\ue_k^\bot=(-u_k^2,u_k^1)$ (  the function  $\gamma_k$ was introduced by  Amick in the paper~\cite{Amick}\,).
By direct calculation,
$$\nabla\gamma_k=-\nabla^\bot\omega_k-\psi_k\nabla\omega_k.$$
Then
\begin{equation}
\label{bpr8-}
\nabla\gamma_k\cdot\nabla^\bot\omega_k=-|\nabla\omega_k|^2.
\end{equation}
In other words,
\begin{equation}
\label{exp-11}
\frac{\partial\gamma_k}{\partial s}:=\nabla\gamma_k\cdot\frac{\nabla^\bot\omega_k}{|\nabla\omega_k|}\equiv-|\nabla\omega_k|,
\end{equation}
where we denote by $\frac{\partial\gamma_k}{\partial s}$ the derivative of~$\gamma_k$ with respect to the~direction tangent to the level set~$\omega_k=c$.
The last identities imply  the~following   monotonicity properties
\begin{equation}
\label{exp-app2}
\begin{array}{lcr}
 \;\;\;\;\;\mbox{$\gamma_k$ \textsl{is monotone along
 level sets of the~vorticity}~$\omega_k=c$ \textsl{and} }\\
\mbox{\textsl{vice versa~-- the~vorticity} $\omega_k$ \textsl{is monotone along level sets of}~$\gamma_k=c$}
\end{array}
\end{equation}
(see \cite{Amick}).  Moreover, there holds the evident identity
\begin{equation}
\label{exp-9}
\gamma_k=\Phi_k\qquad\mbox{  whenever }\ \ \omega_k=0.
\end{equation}

 Suppose~(\ref{bpr8}) is not true.  Let $V_k$ be a connected component of the open set
$\{x\in\Omega^s_k:\omega_k(x)>0\}$. By our assumption, $V_k\ne\Omega^s_k$, therefore,
$\Omega^s_k\cap\partial V_k\ne\emptyset$.
Take a decreasing sequence of noncritical values $\tau_m>0$  of the vorticity~$\omega_k$ satisfying
$\tau_m\to0$ as $m\to\infty$. Since $\omega_k$ is a real analytical function in~$\Omega_k$, for every $\tau_m$ the set
$$\{z\in\partial\Omega^s_k:\omega_k(z)=\tau_m\}$$
is finite. Using this fact and regularity of the values~$\tau_m$, it is easy to see that the set
$$\{z\in\overline{\Omega}^s_k:\omega_k(z)=\tau_m\}$$
is a finite disjoint union of smooth curves homeomorphic to the unit interval~$[0,1]$ with endpoints on~$\partial\Omega^s_k$ (note, that $\tau_m$-level set can not contain curves
 homeomorphic to the circle because of the~monotonicity property~\eqref{exp-app2}\,).

Fix $z_0\in V_k$ and denote by $V_{k,m}$ the sequence of the~connected components of the open set
$\{x\in V_k:\omega_k(x)>\tau_m\}$ containing~$z_0$.
Evidently,
\begin{equation}
\label{en2}
V_{k,m}\subset V_k
\end{equation}
and
\begin{equation}
\label{en3}
V_k=\bigcup\limits_{m\in\N}V_{k,m}
\end{equation}
Now we have to consider three possible cases:

\begin{itemize}
\item[($i$)] \,the equality
\begin{equation}
\label{en4}
S_k(t_1)\cap\partial{V_{k,m}}=\emptyset
\end{equation}
holds for all $m\in\N$;

\item[($ii$)] \,the equality
\begin{equation}
\label{en5}
S_k(t_2)\cap\partial{V_{k,m}}=\emptyset
\end{equation}
holds for all $m\in\N$;

\item[($iii$)] \,the   relations
\begin{equation}
\label{en6}
S_k(t_1)\cap\partial{V_{k,m}}\ne\emptyset\quad\mbox{and} \quad S_k(t_2)\cap\partial{V_{k,m}}\ne\emptyset\end{equation}
hold for all sufficiently large~$m\ge m_0$.
\end{itemize}

Consider the case~(i). First of all, we claim that for this case the relation
\begin{equation}
\label{exp-7}
 S_k(t_2)\cap\partial{V_{k,m}}\ne\emptyset\end{equation}
holds for all
$m$. Indeed, if $S_k(t_2)\cap\partial{V_{k,m}}=\emptyset$, then from~(\ref{en4}) we have
$\partial V_{k,m}\cap\bigl(\partial \Omega^s_k\bigr)=\emptyset$, but this contradicts the~strong  maximum principle for the vorticity~$\omega_k$ (see, $e.g.$, \cite{GW2}\,).

By (\ref{en4}), \eqref{exp-7} there exists evidently a~smooth   arc $L_m\subset V_k\cap\partial V_{k,m}$ such that
\begin{equation}
\label{en7}
\omega_k\equiv \tau_m\qquad\mbox{ \textsl{on} }L_m,\end{equation}
\begin{equation}
\label{en8}
L_m\mbox{ \ \textsl{separates the point} $z_0$ \textsl{from the cycle}~$S_k(t_1)$,\ }\end{equation}
\begin{equation}
\label{en90}
\mbox{ \ \textsl{the endpoints} $A_m$, $B_m$ \textsl{of the arc} $L_m$ \textsl{belong to the cycle}~$S_k(t_2)$}.
\end{equation}
(By {\it cycle} we mean the plane curve which is homeomorphic to the unit circle.) In particular, from the assertion~(\ref{en8})  it follows, that
\begin{equation}
\label{en9}
\mbox{ \ \textsl{the arc}~$L_m$ \textsl{does not degenerate to a point when}~$m\to\infty$}.
\end{equation}
In other words,
\begin{equation}
\label{en9---}
\diam L_m\nrightarrow 0\mbox{ \ \ \ as }m\to\infty.
\end{equation}
By construction we have
$$\gamma(B_m)-\gamma(A_m)=\Phi_k(B_m)-\Phi_k(A_m)-\tau_m(\psi(B_m)-\psi(A_m))
$$
$$
=-\tau_m(\psi(B_m)-\psi(A_m))\to0
$$
as $m\to\infty$. From this fact and  from the~properties~(\ref{en9---}), (\ref{exp-11})--(\ref{exp-app2})  (which could be applied because of (\ref{en7})\,) it follows that
$$\sup\limits_{z\in L_m}|\nabla\omega_k(z)|\to 0\qquad\mbox{ as }m\to\infty.$$
But the last assertion, in view of~(\ref{en9})--(\ref{en9---}), contradicts the fact that $\omega_k$ is a nonconstant real analytical function (in particular, $\nabla\omega_k$ can not  be identically
zero on a~compact connected set which is not a single point).

The case~(ii) can be  proved exactly by the same arguments.

Consider  the last possible case~(iii).
In this case  evidently there exist two smooth arcs $L_m^+$ and $L_m^-$ with the following properties

\begin{itemize}
\item[($\circ$)] \,the arcs $L_m^+$ and $L_m^-$ are homeomorphic to the closed unit interval~$[0,1]$,   their   endpoints $A_m^+,B_m^+$ and $A_m^-,B_m^-$  satisfy
\begin{equation}
\label{en10} A_m^+, A_m^-\in S_k(t_1); \qquad B_m^+,B_m^-\in S_k(t_2).
\end{equation}

\item[($\circ\circ$)] \,the identity
\begin{equation}
\label{en11} \omega_m(z)\equiv \tau_m
\end{equation}
holds  for all $z\in L_m^+$ and $L_m^-$;

\item[($\circ\circ\circ$)] \,the function $\gamma$ is increasing along $L_m^+$ in the direction from~$A_m^+$ to $B_m^+$;
and  $\gamma$  is decreasing along $L_m^-$ in the direction from~$A_m^-$ to $B_m^-$.
\end{itemize}

From the last property it follows that
$$0\le \gamma(A^-_m)-\gamma(B^-_m)=\Phi_k(A^-_m)-\Phi_k(B^-_m)-\tau_m(\psi(A^-_m)-\psi(B^-_m))
$$
$$
=(t_1-t_2)-\tau_m(\psi(B_m^-)-\psi(A_m^-)).
$$
Since by construction $t_2>t_1$ and $\tau_m\to0$ as $m\to\infty$,  we   conclude  that the right hand side of the last   formula is strictly negative for sufficiently
large~$m$, and that is a~required contradiction.

  Thus, the property~(\ref{bpr8}) is proved. Without loss of generality we  may assume that
\begin{equation}
\label{en14}
\omega_k(z)>0\qquad\forall z\in\Omega^s_k.
\end{equation}

\

{\sc Step 7.}  For an angle $\theta\in(0,2\pi)$ denote by $L_\theta$ the ray starting from the origin:
$$L_\theta=\{s(\cos\theta, \sin\theta):\ s\in\R_+\}.$$ Because of assumptions on smallness of the  integrals~(\ref{nont-1}), (\ref{nont-5}) there exists a~value~$\tl$ such that
\begin{equation}
\label{en99}
|\tl-\frac32\pi|<\frac19,
\end{equation}
 and
\begin{equation}
\label{en16}
\int\limits_{L_{\tl}\cap\Omega_k}r|\nabla\u_k|^2(r,\tl)\,dr\le\e_k,
\end{equation}
\begin{equation}
\label{en17}
\int\limits_{L_{\tl}\cap\Omega_k}r|\nabla\tilde\u_k|^2(r,\tl)\,dr\le\e_k,
\end{equation}
\begin{equation}
\label{en18}
\int\limits_{L_{\tl}\cap\Omega'_k}r^2|\nabla\omega_k|^2(r,\tl)\,dr\le\e_k,
\end{equation}
where we denote $\tilde\u_k(r,\theta):=\u_k(r,\theta)-\bar\u_k(r)$ and as usual $\e_k\to0$.
By estimate~(\ref{nont-19}),
$$\sup\limits_{m=1,\dots, M}|\tilde\u_k(R_{mk},\tl)|\le\e_k.$$
Therefore,  the inequality~(\ref{en17}) and   the assumption~(\ref{nont-17})
  yield
\begin{equation}
\label{en19}
\sup\limits_{R_{m1}\le r\le R_{Mk}}|\tilde\u_k(r,\tl)|\le\e_k.
\end{equation}
Then from the identities~(\ref{ess-angle-3})--(\ref{ess-angle-4})
we conclude, that
\begin{equation}
\label{en20}
\u_k(r,\tl)=f_k(r)(\cos\varphi_k(r),\sin\varphi_k(r)),
\end{equation}
where
\begin{equation}
\label{en21}
f_k(r)=|\u_k(r,\tl)|>0,\qquad \qquad|\varphi_k(r)| = \e_k\qquad\forall r\in\bigl[ R_{m1},R_{Mk}\bigr].
\end{equation}

\

{\sc Step 8.}  By construction, the considered domain $\Omega^s_k$ is contained in the   ring
\begin{equation}
\label{en15}
R_{1k}\le |z|\le R_{Mk}\qquad\forall z\in\Omega^s_k
\end{equation}
(see~(\ref{lsa})\,).

Consider the intersection of the~annulus--shape domain $\Omega^s_k$ with the ray~$L_\tl$.
Take a segment $[A,B]\subset L_\tl$ with  endpoints $A,B$
satisfying
\begin{equation}
\label{en22}
|A|=r_1<r_2=|B|,
\end{equation}
\begin{equation}
\label{en23}
A\in S_k(t_1),\qquad B\in S_k(t_2),
\end{equation}
\begin{equation}
\label{en24}
\mbox{\textsl{the interior of the segment}~$[A,B]$ \textsl{is contained in the set}~$\Omega^s_k\cap L_\tl$}
\end{equation}
(the existence of such segment is geometrically evident since $\partial \Omega^s_k=S_k(t_1)\cup S_k(t_2)$ and
the~cycle $S_k(t_2)$ surrounds the cycle~$S_k(t_1)$\,).

Then, by construction and  \eqref{en24-int}, we have
\begin{equation}
\label{en25}
\begin{array}{c}
0<t_2-t_1=\Phi_k(B)-\Phi_k(A)=\int\limits_{[A,B]}\nabla\Phi_k\cdot\et\,dr\\[12pt]
=-\int\limits_{[A,B]}\nabla^\bot\omega_k\cdot\et\,dr+
\int\limits_{[A,B]}\omega_k\u_k^\bot\cdot\et\,dr=I+II,
\end{array}
\end{equation}
where ~$\et=(\cos\tl,\sin\tl)$.
Estimate the~terms~$I$ and $II$ separately:
\begin{equation}
\label{en26}
I\le\sqrt{\int\limits_{[A,B]}r^2|\nabla\omega_k|^2\,dr}\ \ \sqrt{\int\limits_{[A,B]}\frac1{r^2}\,dr}\overset{(\ref{en18})}<
\e_k.
\end{equation}
By formulas~(\ref{en20})--(\ref{en21}), we have
\begin{equation}
\label{en27}
\u_k^\bot(r,\tl)=f_k(r)(\cos\tvv_k(r),\sin\tvv_k(r)),
\end{equation}
where
\begin{equation}
\label{en28}
\tvv_k(r)=\frac\pi2-\varphi_k(r).\end{equation}
Consequently,
\begin{equation}
\label{en29}
\frac\pi2-\e_k<\tvv_k(r)<\frac\pi2+\e_k.\end{equation}
By construction (see~(\ref{en99})\,),
we  also  obtain
\begin{equation}
\label{en30}
\pi-\frac19-\e_k<\tl-\tvv_k(r)<\pi+\frac19+\e_k.
\end{equation}
Finally, we   conclude
\begin{equation}
\label{en31}
\u_k^\bot(\tl,r)\cdot\et<0\qquad\forall r\in \bigl(R_{m1},R_{Mk}\bigr)\supset\bigl[|A|,|B|\bigr].\end{equation}
Therefore, from~(\ref{en14})   it follows that
\begin{equation}\label{en39}
\omega_k\u_k^\bot(\tl,r)\cdot\et<0\qquad\forall r\in\bigl[|A|,|B|\bigr].\end{equation}
Consequently,  the second term~$II$ in ~(\ref{en25}) is negative:
$$II<0,$$
 and ~(\ref{en25})--(\ref{en26}) imply the inequality
$$t_2-t_1<\e_k,$$
 contradicting the choice of~$t_2,t_1$ (recall,  that
$t_2-t_1\ge \frac14-\frac1{40}=\frac9{40}>\frac15,$
see Step~5). This contradiction finishes the proof of~Theorem~\ref{TT2}.
\qed

\

The proof of Theorem~\ref{TT3}   is analogous.
The only differences are   the~following.  The~pressure again is "almost zero", i.e.,  $|p_k|= \e_k$ in the subdomain
$\Omega'_k=B_{\frac34 R_k}$, but
$\u_k\equiv(1,0)$ on $\partial\Omega$ \ and \  $\u_k$ is   zero on the big
circle $S_{R_k}$. Denote
$$R_{0k}=\max\biggl\{r\le R_k:\biggl|\dashint_{S_r}\u_k\,ds\biggr|=\frac15\biggr\}.$$
By the same reasons  as above,
\begin{equation}\label{nont-9---}
R_{0k}\to+\infty\qquad\mbox{ and }\qquad \frac{R_{k}}{R_{0k}}\to+\infty
\end{equation}
as $k\to\infty$.  Then by the maximum principle for the~Bernoulli pressure,
\begin{equation}\label{nont-15----}
\mbox{ \textsl{the function} \ \ }[R_0,R_{0k}]\ni R\mapsto \max\limits_{x\in S_{R}}\Phi_k(x)\quad\mbox{\textsl{is strictly decreasing}},
\end{equation}
where $R_0$ is some fixed radius with $B_{R_0}\supset\partial\Omega$  \,($R_0$ does not depend on~$k$). Then Steps 3--8 of the proof   repeat  exactly   the~corresponding steps    in
the~proof of Theorem~\ref{TT2} with the following obvious change: now the Bernoulli pressure~$\Phi_k$ is decreasing with respect to~$R$, so  that the~circle--type level set $S_k(t_1)\subset\{\Phi_k=t_1\}$  surrounds the level set~$S_k(t_2)
\subset\{\Phi_k=t_2\}$ for $t_1<t_2$.

\

\end{document}